\documentclass{article}

\usepackage{amsmath}
\usepackage{amssymb}

\long\def\rests#1{}

\def\noi{\noindent}

\def\Pf{\noi{\bf Proof.\ \,}}
\def\eop{\hfill\framebox[2.4mm][t1]{\phantom{x}} \vskip 0.15cm } 

\def\voa{vertex operator algebra\ }
\def\gkm{generalized Kac-Moody algebra\ }

\def\a{\alpha}
\def\s{{\bf s}}
\def\g{{\bf g}}
\def\k{{\bf k}}

\def\II{I\!I_{1,1}}
\def\VII{V_{\II}}
\def\Z{{\bf Z}}
\def\z{{\bf z}}

\def\d2{D_{12}(2)}

\def\C{{\cal C}}
\def\D{{\cal D}}
\def\Z{{\bf Z}}
\def\Q{{\bf Q}}
\def\R{{\bf R}}
\def\al{{\alpha}}
\def\Sym{{\rm Sym}}
\def\Aut{{\rm Aut}}

\newcommand{\sfr}[2]{\leavevmode\kern-.05em
  \raise.5ex\hbox{\the\scriptfont0 #1}\kern-.1em
  /\kern-.15em\lower.25ex\hbox{\the\scriptfont0 #2}\kern.02em}
\def\nhalf{{\textstyle \sfr{1}{2}}}
\def\nquart{{\textstyle \sfr{1}{4}}}
\def\ndreiquart{{\textstyle \sfr{3}{4}}}

\title{A generalized Kac-Moody algebra \linebreak of rank~$14$}

\author{Gerald~H\"ohn\thanks{
  Department of Mathematics, Kansas State University, 
  138 Cardwell Hall, Manhattan, KS 66506-2602, USA.
   E-mail: {\tt gerald@math.ksu.edu}}
\ and Nils~R.~Scheithauer\thanks{ 
Fachbereich Mathematik, Technische Universit\"at Darmstadt,
Schlo{\ss}gartenstra{\ss}e 7, 64289 Darmstadt, Deutschland.
E-mail: {\tt scheithauer@mathematik.tu-darmstadt.de} 
\newline
The first author acknowledges financial support from
the Kansas NSF EPSCoR grant NSF32239KAN32240.}
}

\date{}

\begin{document}

\bibliographystyle{amsalpha} 

\newtheorem{thm}{Theorem}[section]
\newtheorem{prop}[thm]{Proposition}
\newtheorem{lem}[thm]{Lemma}
\newtheorem{rem}[thm]{Remark}
\newtheorem{cor}[thm]{Corollary}
\newtheorem{conj}[thm]{Conjecture}
\newtheorem{de}[thm]{Definition}
\newtheorem{nota}[thm]{Notation}

\maketitle

\begin{abstract}
\noindent
We construct a vertex algebra of central charge~$26$
from a lattice orbi\-fold vertex operator algebra of central charge~$12$. 
The BRST-cohomology group of this vertex algebra is a 
new generalized Kac-Moody algebra of rank~$14$. 
We determine its root space multiplicities and a set of simple roots.
\end{abstract}


\section{Introduction}

Generalized Kac-Moody algebras are natural generalizations of the finite dimensional simple Lie algebras 
defined by generators and relations. The denominator identities of these Lie algebras are sometimes 
automorphic forms on orthogonal groups. In this case partial classification results have been 
obtained~\cite{Gritsenko-Nikulin,Scheithauer-classification}.
It is believed that the generalized Kac-Moody algebras whose denominator identities are automorphic 
products of singular weight~\cite{Bo-theta} can be realized as bosonic strings moving on suitable orbifolds. 
The present paper adds evidence to this conjecture. 
Generalized Kac-Moody algebras with such realizations 
can be used to study classification questions in the theory of vertex operator algebras. 
The most prominent application of generalized Kac-Moody algebras 
is Borcherds' proof of the monstrous moonshine conjecture~\cite{Bo-lie}.

So far, there are four generalized Kac-Moody algebras for which explicit
vertex operator algebra constructions are known and the simple roots are determined.
Besides the fake monster Lie algebra~\cite{Bo-fake} and monster
Lie algebra~\cite{Bo-lie} constructed by Borcherds, these
are the fake baby monster Lie algebra constructed by the authors~\cite{HS-fakebaby}
and the baby monster Lie algebra constructed by the first author~\cite{Ho-babymoon}.
The examples studied in~\cite{CKS-lie} depend on the existence of certain vertex operator 
algebras from~\cite{ANS}. The general orbifold approach in~\cite{Carnahan} 
uses hypotheses which are currently unproven. 

The fake monster Lie algebra has rank~$26$ and the fake baby monster Lie algebra has rank~$18$.
{}From several considerations~\cite{Scheithauer-classification, Scheit-conway1,Barnard-theta,
ANS}, 
we expect that the next rank for which generalized Kac-Moody algebras with a natural 
vertex operator algebra construction exist is~$14$. 
In that case, we believe two such algebras exist: One is a ${\bf Z}_3$-twist 
of the fake monster Lie algebra 
which belongs to a series of generalized Kac-Moody algebras investigated 
in~\cite{Bo-lie,Niemann,Scheit-conway1, CKS-lie}. 
The other can be obtained from 
a $\Z_2$-twist of the fake monster Lie algebra corresponding to a class of involutions in 
the isomorphism group of the Leech lattice with a $12$-dimensional fixed point lattice.
In this note, we give a vertex operator algebra construction of this 
generalized Kac-Moody algebra and determine its simple roots. 
The approach of this paper is similar to the one in~\cite{HS-fakebaby}.

\smallskip

This new generalized Kac-Moody algebra together with the fake monster Lie algebra, the fake baby monster
Lie algebra and the monster Lie algebra are the only generalized Kac-Moody algebras 
which can be obtained from a vertex operator algebra associated to a Niemeier lattice 
or the standard $\Z_2$-twist of such a vertex operator algebra \cite{DGM}.
Furthermore, these seem to be all the generalized Kac-Moody algebras  
which can be described by framed vertex operator algebras. 

\medskip

The first three generalized Kac-Moody algebras mentioned in the second paragraph are obtained 
in the following way: Let $V$ be the \voa (VOA) $V_{\Lambda}$ 
associated to the Leech lattice $\Lambda$, the moonshine module VOA 
$V^{\natural}$ or the $\Z_2$-twist of $V_K$,
where $K$ is the Niemeier lattice with root lattice~$A_3^8$.
Let $\VII$ be the vertex algebra of the
two-dimensional even unimodular Lorentzian lattice $\II$. The tensor product
$V\otimes \VII$ is a vertex algebra of central charge $26$.
By using the bosonic ghost vertex superalgebra 
$V_{\rm ghost}$ of central charge $-26$ one defines the
Lie algebra $\g$ as the BRST-cohomology group 
$H^1_{\rm BRST}(V\otimes \VII)$ (cf.~\cite{FGZ}).

In the construction of our new Lie algebra~$\g$,
we take for $V$ a VOA of central charge~$24$ which is obtained by gluing together
the lattice VOA for the rescaled root lattice $\sqrt{2} D_{12}$ 
with the lattice $\Z_2$-orbifold $V_K^+$ of the lattice $K=\sqrt{2}D_{12}^+$.
The decomposition of $V$ into $V_{\sqrt{2}D_{12}}$-modules can be 
described combinatorially using the theory of lattice $\Z_2$-orbifolds
as developed in~\cite{AD-modules, ADL-fusion, Shimakura1}.
This combinatorial description together with the no-ghost
theorem from string theory gives the root lattice and root multiplicities
of~$\g$. Then we construct an automorphic form
on the Grassmannian of negative definite $2$-planes in~$\R^{14,2}$ 
using Borcherds' singular theta correspondence.
The automorphic product can be interpreted as one side of the denominator identity of~$\g$. 
This allows us to determine the simple roots. 

One property which distinguishes $\g$ from the other three examples is that the Weyl group 
of the Lie algebra is not the full reflection group of the root lattice.

\medskip

The paper is organized as follows:
In Section~\ref{voav}, the construction of the \voa $V$ is described
and the $V_{\sqrt{2}D_{12}}$-module decomposition is used to express the
${U}(1)^{12}$-equivariant character of $V$ through theta series 
of the lattice $\sqrt{2}D_{12}$ and a vector valued modular form of weight $-6$.
In the last section, the root lattice, the root multiplicities and 
the simple roots of $\g$ are determined. 

\medskip

We thank Yi-Zhi Huang for helpful discussions.


\section{The \voa $V$ of central \protect\linebreak charge~$24$} 
\label{voav}

In this section, we define a \voa $V$ of central charge~$24$ by gluing together
the lattice \voa~$V_{\sqrt{2}D_{12}}$ 
with the $\Z_2$-orbifold \voa $V_{K}^+$ where $K=\sqrt{2}D_{12}^+$. 
Then we compute its character as a representation for the 
natural Heisenberg subalgebra of~$V$.


\subsection{\bf The \voa $V_N$ and its intertwining algebra}\label{vn}

Let $L\subset {\bf R}^n$ be an even lattice of rank $n$ and let 
$L'=\{\lambda \in {\bf R}^n \, | \, (\lambda,\mu) \in {\bf Z}$ 
\linebreak for all $\mu\in L  \}$ be the dual lattice. 
The map $q_L: L'/L\to \Q/\Z$, 
$\lambda + L \mapsto \linebreak (\lambda,\lambda)/2 \pmod{\Z}$,
gives the discriminant group $L'/L$ the structure of a finite quadratic space which
is called the discriminant form of $L$.
We sometimes write $\lambda^2=(\lambda,\lambda)$ for the norm of $\lambda$.

The isomorphism classes of irreducible modules of the
vertex operator algebra $V_L$ associated to an even 
integral lattice $L$ can be parameterized by the discriminant group $L'/L$
of the lattice~\cite{DoLe}. 
For each coset $\lambda + L  \in L'/L$ there exists a unique 
irreducible $V_{L}$-module which we denote by $V_{\lambda + L}$.

The fusion product between the irreducible modules is given by
$$  V_{\lambda + L} \times  V_{\mu + L} = V_{\lambda + \mu + L} $$
with $\lambda + L$, $\mu + L$ in $L'/L$, i.e. the
fusion algebra of $V_L$ is isomorphic to the group ring ${\bf C}[L'/L]$
and each simple module is a simple current.

The direct sum of the irreducible modules of a lattice vertex operator algebra $V_L$ 
admits the structure of an abelian intertwining algebra~\cite{DoLe}, Th.~12.24, such that
the cohomology class of the associated $3$-cocycle is determined by the quadratic form $q_L$ 
on $L'/L$. The conformal weights modulo~$\Z$ of the irreducible $V_L$-modules $V_{\lambda + L}$,
$\lambda + L \in L'/L$ are the values of the quadratic form $q_L$.

We collect these results in the following theorem.
\begin{thm}\label{abelian-vn}
The direct sum of the simple modules of $V_{L}$ has the structure of an abelian intertwining algebra. 
The associated quadratic space can be identified with the discriminant form $L'/L$.
\eop
\end{thm}

For the proof of some identities, it is useful to
interpret an element $f$ in ${\bf C}[L][[q^{1/k}]][q^{-1/k}]$,
where $L$ is a lattice and $k\in{\bf N}$, 
as a function on ${\cal H}\times (L\otimes {\bf C})$,
where ${\cal H}=\{z\in {\bf C} \mid \text{Im}(z)>0\}$ is the complex
upper half plane. This is done by the substitutions 
$q\mapsto e^{2\pi i \tau}$ and $e^{\s} \mapsto e^{2\pi i(\s,\z)}$
for $(\tau,\z)\in {\cal H}\times (L\otimes {\bf C})$ 
(in the case of convergence). We indicate this by writing $f(\tau,\z)$.

Let $\eta(\tau)=q^{1/24}\prod_{k=1}^{\infty}(1-q^k)$ 
be the Dedekind eta function.
We define 
the theta function of the coset $\lambda + L$ by
$$\theta_{\lambda + L}= \sum_{{\bf s}\in \lambda + L} q^{\s^2/2}\,e^{\bf s}.$$

The ${\bf Z}$-grading on a VOA $W=\bigoplus_{k=0}^{\infty}W_k$ is given 
by the eigenvalues of the Virasoro generator~$L_0$.
Suppose there is an action of a connected compact Lie group $G$ on $W$ respecting this grading. 
Let $L$ be the weight lattice of a maximal torus of $G$.
Then we denote by $W_k(\s)$ the subspace of $W_k$ of weight $\s$.
The character of~$W$ is defined by 
$$\chi_W=q^{-c/24}\sum_{k\in \bf Z} \, \sum_{\s\in L} \dim W_k(\s)\,q^k\,e^{\s},$$
where $c$ is the central charge of $W$.

On $V_L$ there is the action of ${\bf R}^{n}/L'$ by \voa automorphisms.
The $|L'/L|$-fold cover $T={\bf R}^{n}/L$ acts also on the modules $V_{\lambda + L}$
and the weights form the coset $\lambda + L$.

{}From the construction of $V_{\lambda + L}$ one obtains the following
description of the $T$-equivariant graded character.

\begin{lem}\label{gitterchar}
The  $V_{L}$-module $V_{\lambda + L}$ has the character 
$\theta_{\lambda + L}(\tau,{\bf z})/\eta(\tau)^n$. \eop
\end{lem}

\medskip

We choose now for $L$ the lattice $N= \sqrt{2}D_{12}$, i.e.
\[ N = \bigl\{ \sqrt{2}\,(x_1,\,\ldots,\,x_{12})  \in \R^{12} \mid 
\text{all $x_i\in \Z$ and $\sum_{i=1}^{12} x_i \equiv 0 \!\!  \pmod 2$} \bigr\} \, .  
\]
Then the automorphism group of $N$ is generated by the permutations of the coordinates 
and arbitrary sign changes. It has shape $2^{12}.\Sym_{12}$.

\begin{lem}\label{quad-N}
The discriminant group of $N$ and the orbits under the induced action of 
${\rm Aut}(N)\cong 2^{12}.\Sym_{12}$ on the discriminant group are as described in 
Table~\ref{N-discriminant}.
The lattice $N$ has genus $I\!I_{12,0}(2_{I\!I}^{-10} 4_{I\!I}^{-2})$
in the notation of~\cite{CoSl}.
\end{lem}

\Pf
The dual lattice of $N$ is given by
\[ N' = \bigl\{ \tfrac{1}{\sqrt{2}}\,(x_1,\,\ldots,\,x_{12}) \in \R^{12} \mid
\text{all $x_i\in \Z$ or all $x_i\in \Z + \tfrac{1}{2}$} \bigr\}  .
\]
It is easy to describe the decomposition of the discriminant group $N'/N$ 
into orbits of ${\rm Aut}(N)$ and to determine representatives.
The genus can be calculated by diagonalizing a Gram matrix of $N$ over the $2$-adic integers.
Note that the genus is uniquely determined by $(N'/N,q_N)$ and the signature of $N$.
\eop

\begin{table}\caption{Orbits of ${\rm Aut}(N)$ on the discriminant form of $N$.}\label{N-discriminant}
$$\begin{array}{r|ccc|rcc}
\hbox{No.} & \hbox{representative} &  \hbox{$N'$-orbit size} &  \hbox{norm} &  
\hbox{orbit size} &  \hbox{$q_N$} &  \hbox{order} \\ \hline
1 & \frac{1}{\sqrt{2}}\, (0^{12}) & 1 & 0 & 1 & 0 & 1 \\[1mm] \hline
2 & \frac{1}{\sqrt{2}}\,(2,0^{11})  & 2\cdot 12 & 2 & 1 & 0 & 2  \\
3 & \frac{1}{\sqrt{2}}\,(1^{12}) & 2^{12} & 6 & 2 & 0 & 2  \\[1mm] \hline
4 & \frac{1}{\sqrt{2}}\,(1^{4},0^8)  & 2^4{ 12 \choose 4} & 2 & 990 & 0 & 2  \\
5 & \frac{1}{\sqrt{2}}\,(1^{8},0^4) & 2^8{ 12 \choose 8} & 4 & 990 & 0 & 2  \\ [1mm]\hline
6 & \frac{1}{\sqrt{2}}\,(1^{2},0^{10}) & 2^2{ 12 \choose 2} & 1 & 132 & \nhalf & 2  \\
7 & \frac{1}{\sqrt{2}}\,(1^{6},0^6) & 2^6{ 12 \choose 6} & 3 & 1848 & \nhalf & 2  \\
8 & \frac{1}{\sqrt{2}}\,(1^{10},0^2) & 2^{10}{ 12 \choose 10} & 5 & 132 & \nhalf & 2  \\[1mm] \hline
9 & \frac{1}{\sqrt{2}}\,(1,0^{11}) & 2\cdot 12 & 1/2 & 24 & \nquart & 4 \\
10 & \frac{1}{\sqrt{2}}\,(1^{5},0^{7}) & 2^5 { 12 \choose 5}  & 5/2 & 1584 & \nquart & 4 \\
11 & \frac{1}{\sqrt{2}}\,(1^{9},0^{3}) & 2^9{ 12 \choose 9}  & 9/2 & 440 & \nquart & 4 \\ 
12 & \frac{1}{\sqrt{2}}\,(\frac{3}{2},(\frac{1}{2})^{11}) &  2^{12}\cdot  12   & 5/2 & 4096 & \nquart & 4 \\[1mm] \hline
13 & \frac{1}{\sqrt{2}}\,(1^{3},0^{9}) & 2^3 { 12 \choose 3} & 3/2 & 440 & \ndreiquart & 4 \\
14 & \frac{1}{\sqrt{2}}\,(1^{7},0^{5}) & 2^7 { 12 \choose 7}  & 7/2 & 1584 & \ndreiquart & 4 \\
15 & \frac{1}{\sqrt{2}}\,(1^{11},0) & 2^{11} { 12 \choose 11}  & 11/2 & 24 & \ndreiquart & 4 \\ 
16 & \frac{1}{\sqrt{2}}\,( (\frac{1}{2})^{12} )  & 2^{12}  & 3/2 & 4096 & \ndreiquart & 4 \\ [1mm] \hline
\end{array}$$
\end{table}

The first eight orbits of Table~\ref{N-discriminant} form the $2$-torsion subgroup of $N'/N$
 and the first three orbits
consist of elements which are a multiple of $2$ of another element. Thus the orbits in 
Table~\ref{N-discriminant} separated by horizontal lines belong also to different orbits
of $N'/N$ under the action of the automorphism group of the discriminant form of $N$.


\subsection{\bf The \voa $V_K^+$ and its intertwining algebra}\label{vm}

As in the previous subsection, let $L\subset {\bf R}^n$  be an even lattice of rank $n$. We are
interested in the VOA $V_L^+$, the fixed point subspace of $V_L$ under the involution
induced from the $-1$ isomorphism of $L$.
The irreducible modules of $V_L^+$ have been described in~\cite{AD-modules}, their fusion rules in~\cite{ADL-fusion} and
the automorphism group in~\cite{Shimakura1,Shimakura2}. 

We specialize the discussion here to the case of the lattice $K=\sqrt{2}D_{12}^+$, i.e.
\[ K = \bigl\{ \sqrt{2}\,(x_1,\,\ldots,\,x_{12}) 
\mid \text{all $x_i\in \Z$ or all $x_i\in \Z + \tfrac{1}{2}$ 
and $\sum_{i=1}^{12} x_i \equiv 0\!\! \pmod 2$} \bigr\}  .  \]
The automorphism group of $K$ is isomorphic to the Weyl group $W(D_{12})$.
It is generated by permutations of the coordinates and even sign changes  
and has shape $2^{11}.\Sym_{12}$. 
The lattice $D_{12}^+$ is the unique indecomposable unimodular integral lattice in dimension~$12$.

\smallskip

The following lemma is easy to prove.
\begin{lem} 
The discriminant group of $K$ and the orbits under the induced action of 
${\rm Aut}(K)\cong 2^{11}.\Sym_{12}$ on the discriminant group are as described in 
Table~\ref{K-discriminant}.
The lattice $K$ has genus~$I\!I_{12,0}(2_{4}^{+12})$. \eop
\end{lem}

\begin{table}\caption{Orbits of ${\rm Aut}(K)$ on the discriminant form of $K$.}\label{K-discriminant}
$$\begin{array}{r|ccc|rcc}
\hbox{No.} & \hbox{representative} &  \hbox{$K'$-orbit size} &  \hbox{norm} &  
\hbox{orbit size} &  \hbox{$q_K$} &  \hbox{order} \\[1mm] \hline
1 &\frac{1}{\sqrt{2}}\, (0^{12}) & 1 & 0 & 1 & 0 & 1 \\ [1mm] \hline
2 &\frac{1}{\sqrt{2}}\, (2,0^{11})  & 2\cdot 12 & 2 & 1 & 0 & 2  \\[1mm] \hline
3 &\frac{1}{\sqrt{2}}\,(1^{4},0^8)  & 2^4\cdot {12 \choose 4} & 2 & 990 & 0 & 2  \\[1mm] \hline
4 &\frac{1}{\sqrt{2}}\,(1^{2},0^{10})  & 2^2\cdot {12 \choose 2} & 1 & 132 & \nhalf & 2  \\
5 & \frac{1}{\sqrt{2}}\,(1^{6},0^6) & 2^6\cdot {12 \choose 6} & 3 & 924 & \nhalf & 2  \\[1mm] \hline
6 & \frac{1}{\sqrt{2}}\,(-\frac{3}{2},(\frac{1}{2})^{11}) &  2^{11}\cdot  12   & 5/2 & 1024 & \nquart & 2  \\[1mm] \hline
7 & \frac{1}{\sqrt{2}}\,( (\frac{1}{2})^{12} )  & 2^{11}  & 3/2 & 1024 & \ndreiquart & 2  \\
\end{array}$$
\end{table}

Let $\tau$ be the involution in  ${\rm Aut}(V_K)$ which is up to 
conjugation the unique lift of the involution $-1$ in ${\rm Aut}(K)$ 
to ${\rm Aut}(V_K)$ (cf.~\cite{DGH-virs}, Appendix~D). Denote by  $V_K^+$ the
fixed point vertex operator subalgebra of $V_K$ under the action of $\tau$.

The isomorphism classes of irreducible modules of $V_K^+$ are described in~\cite{AD-modules}. Since
$K$ is $2$-elementary, that is $2K'\subset K$, the discussion can be simplified, cf.~\cite{Shimakura1},
Section~3.2. The isomorphism classes of irreducible modules of $V_K^+$ consist of the 
so called untwisted modules $V_{\lambda + K}^\pm$ where $\lambda + K$ runs
through the discriminant group $K'/K$ and certain so called twisted modules $V_K^{T_\chi,\pm}$.

The fusion rules between the $2^{13}$ modules $V_{\lambda + K}^{\pm}$ are
$$ V_{\lambda + K}^\delta \times V_{\mu+K}^\epsilon = V_{\lambda + \mu + K}^{\pm}, $$
where $\delta$, $\epsilon\in \{\pm\} \cong\Z_2$ and $\lambda$, $\mu\in K'$
and the exact sign in $V_{\lambda + \mu + K}^{\pm}$ can be determined from 
the discriminant form of $K$.

Since the fusion product $\times$ is commutative and associative we see that it induces
on the set $\{ V_{\lambda + K}^\pm \mid  \lambda + K \in K'/K\}$ 
of isomorphism classes of untwisted $V_{K}^+$-modules the structure of an
abelian group of exponent $4$. In fact, $V_{\lambda + K}^\pm$ is of order $4$ if
and only if $\lambda$ has non-integral norm, cf.~\cite{Shimakura1}, Remark~3.5.
More precisely, we have the following description of the fusion algebra of $V_K^+$
(\cite{ADL-fusion}, Theorem 5.1; see also p.~216).
\begin{prop}\label{fusionalgebra}
The isomorphism classes of irreducible modules of $V_K^+$ form an abelian group $A$ of
exponent $4$ under the fusion product which is isomorphic to $\Z_2^{10}\times \Z_4^2$.  \eop
\end{prop}
In particular, all the twisted modules of $V_K^+$ are of order $4$ in $A$.

\smallskip

We need a characterization of modular tensor categories with fusion algebra as 
in Proposition~\ref{fusionalgebra}.
\begin{lem}\label{classabelianmod}
Let ${\cal C}$ be a modular tensor category such that the fusion algebra is isomorphic to the group ring 
of a finite abelian group~$A$. Then the category ${\cal C}$ is up to equivalence determined by 
the twist $\theta$ considered as a map $\theta: A\longrightarrow {\bf C}^*$ such that
$\theta_a=\theta(a)\,{\rm id}\in {\rm Hom}(a,a)$ for $a\in A$. \eop
\end{lem}
A proof can be found in~\cite{FRS04}, Prop.~2.11 (ii)--(iv). Note that the cohomology class of an abelian 
$3$-cocycle $(F,\Omega)$ on $A$ is determined by the corresponding quadratic form \hbox{$q:A\longrightarrow {\bf C}^*$,}
$q(a)=\Omega(a,a)$, (\cite{EM50,MacLane}) and the quadratic form $q$ is the inverse of the twist $\theta$
(cf.~\cite{FRS04}, Prop.~2.14).

The following general result is probably known 
(cf.~\cite{Huang-intertwining}, Prop.~3.4 for a related statement). 
We include a proof for completeness. 
\begin{thm}\label{abintertwiner}
Let $V=\bigoplus_{n=0}^\infty V_n$ be a simple VOA with $V_0 = {\bf C}{\bf 1}$, $V'=V$ and  
for which every weak $V$-module is a direct sum of irreducible $V$-modules. 
Assume that the fusion algebra defined by the fusion rules 
for the isomorphism classes of irreducible $V$-modules is isomorphic to the group ring of a finite abelian
group~$A$. Then the direct sum of representatives $W^a$ of the isomorphism classes of irreducible modules can
be given the structure of an abelian intertwining algebra as in~\cite{DoLe} with associated quadratic space 
$(A,q_A)$ where the quadratic form $q_A$ on $A$ is given by $q_A(a)=h(a)\pmod{\Z}$ with $h(a)$ denoting the 
conformal weight of the module $W^a$.
\end{thm}
\noindent{\bf Proof.} 
It was proven by Huang that a direct sum of representatives of irreducible $V$-modules 
has the structure of an intertwining algebra (see \cite{Huang-DG}, Theorem~3.7).
The associativity and commutativity properties of the intertwining operators allow to define fusing 
and braiding isomorphisms satisfying the pentagon and the two hexagon identities. By choosing a basis 
in each intertwining space they may be represented by matrices.

The assumption on the fusion algebra
guarantees that the intertwining space ${\cal V}_{a_1a_2}^{a_3}$ for three irreducible modules 
$a_1$, $a_2$, $a_3\in A$ is non-zero only if $a_1+a_2=a_3$ in $A$ in which case the space is one-dimensional.
It follows that the non-zero entries for the fusing and braiding matrices can be described by maps
$$F: A\times A \times A\longrightarrow {\bf C}^*\qquad\hbox{and}\qquad B: A\times A \times A\longrightarrow {\bf C}^*.$$
The values of $F$ and $B$ are non-zero since the fusing and braiding maps are isomorphisms.
Define $\Omega:A\times A\longrightarrow {\bf C}^*$ by  $\Omega(a_1,a_2)=B(a_1,a_2,0)$. Then the 
pentagon and the two hexagon identities imply that $(F,\Omega)$ represents a cocycle for the $4$-th cohomology
group of the Eilenberg-MacLane space $K(A,4)$ (see \cite{JS}, Prop.~13 and \cite{DoLe}, Remark~12.22).
A basis change for the intertwining spaces corresponds to adding a coboundary to $(F,\Omega)$ 
(cf.\ \cite{DoLe}, Remark~12.21.). The function $q: A\longrightarrow {\bf C}^* $, $ q(a)=\Omega(a,a)$,
is a quadratic form on $A$. As mentioned above, the cohomology group $H^4(K(A,4),{\bf C}^*)$ is isomorphic to the group
of quadratic forms on $A$ with values in ${\bf C}^*$. We finally note that $(F,\Omega)$
may be assumed to be normalized (\cite{DoLe}, end of Remark~12.22). 

By replacing the Jacobi identity axiom in the definition of an abelian intertwining algebra by the generalized
commutativity and associativity properties (\cite{DoLe}, Theorem~12.32), we see that the properties
of an intertwining algebra as in~\cite{Huang-intertwining} together with the choices of nonzero
intertwining operators in the intertwining spaces (extending the $V$-module structure) imply all the properties 
for an abelian intertwining algebra for the pair $(F,\Omega)$ as in~\cite{DoLe} with the exception of 
Eq.\ (12.120) of \cite{DoLe}. By Remark 12.29 of \cite{DoLe}, the axiom (12.120) can be replaced by the
grading condition $W^a=\bigoplus_{e^{2\pi i n}=q(a)^{-1}} W^a_n$ for $a\in A$.

\smallskip

Huang and Lepowsky defined in \cite{HL-tensorcat} the structure of a braided tensor category on the 
category of $V$-modules (see~\cite{HL-tensorcat}, Theorem~4.4).
The intertwining spaces used above can be identified with the spaces 
${\rm Hom}(W^{a_1}\otimes W^{a_2}, W^{a_3})$ in the tensor category, 
cf.~\cite{HL-tensorcat}, Prop.~3.4.

Under the assumption on the VOA as in the theorem it was shown
in \cite{Huang-modular}, Th.~4.6, that this tensor category has in addition the structure of a
modular tensor category where the twist $\theta_W:W\longrightarrow W$ is given by the operator
$e^{2\pi i L_0}$ (\cite{Huang-modular}, Th.~4.1).


For a ribbon category whose underlying braided tensor category is defined by an abelian group $A$ 
and a quadratic form $q:A\longrightarrow {\bf C}^*$ as above, we obtain from Lemma~\ref{classabelianmod} 
that $\theta_{W^a}=q(a)^{-1}$. 
This is equivalent to the above stated grading condition finishing the proof that $\bigoplus_{a\in A}W^a$ 
has the structure of an abelian intertwining algebra with 
$q_A(a)=-\frac{1}{2\pi i}\log  q(a)=h(a) \pmod{\Z}$.
\eop

\medskip

\begin{thm}\label{abelian-vc}
The direct sum of representatives for the isomorphism classes of irreducible modules of $V_K^+$ can be
given the structure of an abelian intertwining algebra. 
The associated quadratic space $(A,q_A)$ is isomorphic to the discriminant 
form of the lattice $N$.
\end{thm}

{\bf Proof.} 
First we have to check that the vertex operator $V_K^+$ satisfies all the properties assumed 
in Theorem~\ref{abintertwiner}. The grading condition is clear. The simplicity follows from the classification
of simple modules~\cite{AD-modules}. Since $(V_K^+)_1=0$ we have $L_1(V_K^+)_1\not= (V_K^+)_0$
and hence $(V_K^+)'=V_K^+$ by~\cite{li-dual}. In~\cite{abe-rational}, it was proven that a vertex operator algebra of 
type $V_L^+$ with positive definite even lattice $L$ is rational. In~\cite{ABD-regular}, it was shown that such a 
vertex operator algebra $V_L^+$ satisfies the so-called $C_2$-cofiniteness condition. 
It was also shown in~\cite{ABD-regular} 
that a rational vertex operator algebra satisfying the grading and the $C_2$-cofiniteness condition has the property 
that every weak module is the direct sum of ordinary irreducible modules.

We know already from Proposition~\ref{fusionalgebra} that the fusion algebra of $V_K^+$ 
has the required form $A=\Z_2^{10}\times \Z_4^2$. It remains to 
determine the quadratic form $q_A:A\longrightarrow {\bf Q}/{\bf Z}$. 
Since Theorem~\ref{abintertwiner} implies that $q_A$ is determined by the conformal weights of the modules of $V_K^+$,
it is enough to know the characters of those modules. Those will be determined in the following discussion
showing that $q_A$ has indeed the claimed form (see also Table~3).  
\eop

\smallskip 

We consider now the characters of the irreducible $V_K^+$-modules. They will also be used in the next section.
One has~\cite{FLM}
\begin{eqnarray}\label{characters}
 \chi_{V_K^\pm} & = & \frac{1}{2}\left(\frac{\theta_{K}(q)}{\eta(q)^{12}}\pm \frac{\eta(q)^{12}}{\eta(q^2)^{12}}\right),   \nonumber \\
 \chi_{V_{\lambda + K}^\pm} & = &\frac{1}{2}\frac{\theta_{\lambda + K}(q)}{\eta(q)^{12}}, \quad \hbox{for $\lambda + K\not=K$,} \\ \nonumber
 \chi_{V_K^{T_\chi,\pm}} & = &\frac{1}{2}\,q^{3/4}\,\left(\frac{\eta(q)^{12}}{\eta(q^{1/2})^{12}} \pm \frac{\eta(q^2)^{12}\eta(q^{1/2})^{12}}{\eta(q)^{24}}\right) . 
\end{eqnarray}
In particular, the characters of the $V_K^+$-modules $V_{\lambda + K}^{\pm}$ depend for $\lambda + K\not=K$
only on the orbit of $\lambda + K$ under $\Aut(K)$ in $K'/K$. We denote the character of $V_{\lambda + K}^+$ for
$\lambda + K$ belonging to the orbit $n$ in Table~\ref{K-discriminant} by $g_n$.
An explicit computation gives
$$\begin{array}{lcl}
 g_1 & = &q^{-1/2}\,( 1+210 \,q^2+2752 \,q^3+29727 \,q^4+225408 \,q^5+ \cdots),  \\ 
 g_2 & = &q^{-1/2}\,( 12 \,q+144 \,q^2+2984 \,q^3+29088 \,q^4+227004 \,q^5 \cdots),\\
 g_3 & = &q^{-1/2}\,( 4 \,q+176 \,q^2+2872 \,q^3+29408 \,q^4+226196 \,q^5+\cdots), \\
 g_4 & = &q^{-1/2}\,( \,q^{1/2}+32 \,q^{3/2}+768\,q^{5/2}+9600 \,q^{7/2}+83968 \,q^{9/2}+\cdots),\\ 
 g_5 & = &q^{-1/2}\,( 32 \,q^{3/2}+384 \,q^{5/2}+4992\,q^{7/2}+49408 \,q^{9/2}+\cdots),\\ 
 g_6 & = &q^{-1/2}\,( 12 \,q^{5/4}+376 \,q^{9/4}+5316 \,q^{13/4}+50088\,q^{17/4}+\cdots),\\
 g_7 & = &q^{-1/2}\,( \,q^{3/4}+78 \,q^{7/4}+1509 \,q^{11/4}+16966\,q^{15/4}+\cdots).
\end{array}$$

\smallskip 

Now we discuss the automorphism group of $V_K^+$ (see~\cite{Shimakura1}) and its induced action on the
quadratic space $(A,q_A)$ 
although this information is not really necessary for the construction and understanding of the \gkm $\g$.
The centralizer $H$ of $\tau$ in ${\rm Aut}(V_K)$ acts on $V_K^+$. 
$H$ has shape $2^{12}.{\rm Aut}(K)$, where the $2^{12}$ can be identified with
${\rm Hom}(K,\Z_2)$. The element $\tau\in H$ acts trivially. The induced action of $H$ on
the set of isomorphism classes of irreducible $V_K^+$-modules stabilizes the set of
untwisted modules. For $g\in H$ one has
$$g \big( \{V_{\lambda + K}^\pm \} \big)  = \{V_{g(\lambda + K)}^\pm\} \, .$$
Moreover, if $g\in {\rm Hom}(K,\Z_2)\subset H$ then
\[ g \big( V_{\lambda + K}^\pm \big) = 
\begin{cases} 
V_{\lambda + K}^\pm & \text{\ if\ } g(2\lambda)=0, \\
V_{\lambda + K}^\mp & \text{\ if\ } g(2\lambda)=1  \, . 
\end{cases}  \]

Thus if we have an element $\lambda \in K'$ for which $2\lambda$ is not in $2K$, 
i.e.\ $\lambda\not\in K$,
we can find an element $g\in {\rm Hom}(K,\Z_2)\subset H$ with $g(2\lambda)=1$. 
It follows that the modules $V_{\lambda + K}^\pm$, where $\lambda + K$ runs through an ${\rm Aut}(K)$-orbit,
all belong to the same $H$-orbit.

It was shown by Shimakura~\cite{Shimakura1} that the orbit of $V_K^-$ under ${\rm Aut}(V_K^+)$ 
contains in addition to $V_K^-$ the $V_{\lambda + K}^\pm$, $\lambda + K\not=K$,
for which the number of norm $2$ vectors in $\lambda + K$ is exactly $24$ 
(because the dimension of $K$ is different from $8$ or $16$).
This is exactly the case when $\lambda + K$ belongs to the orbit no.~2 of size~$1$ in Table~\ref{K-discriminant}.
{}From this it can be deduced that ${\rm Aut}(V_K^+)$ has shape $2^{11}.2^{10}.\Sym_{12}.\Sym_3$:
Since $K$ can be constructed as $L_C^+$ for the binary code $C=\{0^{12},\,1^{12}\}$
(cf.~Remark~\ref{framedvoa} below), there exists an triality automorphism responsible for the 
$\Sym_3\cong {\rm GL}({\bf F}_2^2)$
and it is enough to know that ${\rm Aut}(C)=\Sym_{12}$. For the details 
cf.~\cite{Shimakura1}, Th.~4.3~(iv) where the case of ${\rm Aut}(V_L^+)$ for
$L\cong{\sqrt{2}D_{12}}$ is discussed.

We collect the results in Table~\ref{VKP-modules}. Here we write $[n]^\pm$ for the set of
modules $V_{\lambda + K}^\pm$ for which $\lambda + K$ belongs to the orbit no.~$n$ in 
Table~\ref{K-discriminant}.
\begin{table}\caption{Orbits of ${\rm Aut}(V_K^+)$ on the irreducible modules of $V_K^+$.}\label{VKP-modules}
$$\begin{array}{r|crcccc}
\hbox{No.} & \hbox{$H$-orbits} & \hbox{orbit size} &  \hbox{\quad $h$\quad } &  \hbox{$q$}  &  \hbox{order} & \hbox{character} \\ \hline
1 & [1]^+ & 1 & 0 & 0 & 1 & g_{1} \\ \hline
2 & [1]^-,\, [2]^+,\, [2]^-    & 3\times 1 & 1 & 0 &2 & g_{2} \\ \hline
3 & [3]^+,\, [3]^-            & 2\times 990 & 1 & 0 & 2 &  g_{3}\\ \hline
4 & [4]^+,\, [4]^-            & 2\times 132 & 1/2 & \nhalf & 2 & g_{4}\\
5 & [5]^+,\, [5]^-         &  2\times 924 & 3/2 & \nhalf & 2 & g_{5}\\ \hline
6 & [6]^+,\, [6]^- ,\,\{[\chi]^- \}    & 3\times2\times 1024 & 5/4 & \nquart & 4 & g_{6} \\ \hline
7 & [7]^+,\, [7]^- ,\,\{[\chi]^+ \}    &3\times2\times 1024 & 3/4  & \ndreiquart & 4 & g_{7} \\ \hline
\end{array}$$
\end{table}
Note that $g_n$ is the character of the $V_K^+$-modules in the $n$-th orbit in Table~\ref{VKP-modules}.

The only entries in Table~\ref{VKP-modules} which remain to be discussed are the $H$-orbits of
the twisted modules. 
If $V_K^+$ is extended by the unique module belonging
to $[2]^+$ or $[2]^-$ then one obtains an extension $\widetilde V$ of the VOA $V_K^+$ which is isomorphic
to the lattice VOA $V_K$, but some twisted $V_K^+$-modules become now untwisted modules 
for $V_K^+$ considered as the fixed point VOA of $\widetilde V \cong V_K$. 
Under the extra automorphisms in ${\rm Aut}(V_K^+)$ which map $[1]^-$ to $[2]^+$ or $[2]^-$, 
a twisted $V_K^+$-module may be mapped to an untwisted one. 
In fact, this can be done for all the twisted $V_K^+$-modules, cf.~\cite{FLM}, Chapter 11. 
Now it follows from the given conformal characters that all twisted modules 
$V_K^{T_\chi,-}$ belong to the orbits $[6]^+$ and $[6]^-$ 
and all twisted modules $V_K^{T_\chi,+}$ belong to the orbits $[7]^+$ and $[7]^-$.

\pagebreak



\begin{rem}\label{orthogonalorbs}
There are exactly $6$ orbits under the action of ${\rm Aut}(A,q_A)$ on $A$. 
\end{rem}

\Pf
We only have to show that the elements of order $2$ and norm $1/2$ in $A$ 
are conjugate under ${\rm Aut}(A,q_A)$, cf.~Table~2 or~3.
Let $\gamma$ be such an element. We choose a Jordan decomposition of $A$. 
Then the projection of $\gamma$ on the summand $2_{I\!I}^{-10}$ is nontrivial 
because $4_{I\!I}^{-2}$ contains no elements of norm $1/2$. The bilinear form is nondegenerate 
on $2_{I\!I}^{-10}$ so that there is an element $\mu$ in $2_{I\!I}^{-10}$ such that $(\gamma,\mu) = 1/2$. 
We can assume that $\mu$ has norm $1/2$. Then $\langle \gamma, \mu \rangle$ is a discriminant 
form of type $2_{I\!I}^{-2}$ and $A = \langle \gamma, \mu \rangle \oplus \langle \gamma, \mu \rangle^{\perp}$. 
The statement now follows from the uniqueness of the $2$-adic symbol of $A$.\eop



\medskip


\begin{rem}\label{framedvoa}
$V_K^+$ is isomorphic to the framed code VOA $V_{\C}$ 
associated to the binary code $\C$ dual to the code $\D$ with generator matrix
$$\left(\begin{array}{c}
1  1  1  1\  1  1  1  1\  0  0  0  0\  0  0  0  0\  0  0  0  0\  0  0  0  0 \\
0  0  0  0\  1  1  1  1\  1  1  1  1\  0  0  0  0\  0  0  0  0\  0  0  0  0 \\
0  0  0  0\  0  0  0  0\  1  1  1  1\  1  1  1  1\  0  0  0  0\  0  0  0  0 \\
0  0  0  0\  0  0  0  0\  0  0  0  0\  1  1  1  1\  1  1  1  1\  0  0  0  0 \\
0  0  0  0\  0  0  0  0\  0  0  0  0\  0  0  0  0\  1  1  1  1\  1  1  1  1 \\
1  1  0  0\  1  1  0  0\  1  1  0  0\  1  1  0  0\  1  1  0  0\  1  1  0  0 \\
1  0  1  0\  1  0  1  0\  1  0  1  0\  1  0  1  0\  1  0  1  0\  1  0  1  0
\end{array}\right)_.$$ 
\end{rem}
\Pf
The lattice $K$ can be written in terms of the binary code $C=\{0^{12},\,1^{12}\}$
of length~$12$ which is generated by the overall-one vector $(1111\,1111\,1111)$ as
$$K=L_C^+=\frac{1}{\sqrt{2}}\{c+x\mid c\in C,\ x\in(2\Z)^{12}\hbox{\ with\ } \sum_{i=1}^{12}x_i \equiv 0 \! \! \pmod 4 \}.$$ 
Now the result follows from the Virasoro decomposition of $\widetilde{V}_{\widetilde{L}_C}$ given 
in~\cite{DGH-virs}, Th.~4.10, by observing that the first term in the sum corresponds to $V_{L^+_C}^+$
so that $V_K^+$ and $V_{\C}$ have the same Virasoro decomposition
and must therefore be isomorphic (see~\cite{DGH-virs}, Prop.~2.16 and~\cite{Ho-genus}, Th.~4.3). 
Note that for $C$ all markings are equivalent and that the
proof of Th.~4.10 in~\cite{DGH-virs} shows that the 
self-duality assumption on $C$ is unnecessary.
\eop

For proving Theorem~\ref{abelian-vc} one can also use this remark and
the results of~\cite{Mi-rep}, where all irreducible modules of a 
framed code VOA $V_{\C}$ are described.


\subsection{\bf The gluing of $V_N$ and $V_K^+$ }

The quadratic spaces $(A,q_A)$ and $(A,-q_A)$ are isomorphic.
We choose an isomorphism $i:N'/N \to A$ between the spaces $(N'/N,q_N)$ and $(A,-q_A)$. 
Let~$V$ be the $V_N\otimes V_K^+$-module
$$V\ =\ \bigoplus_{\lambda \in N'/N}  V_{\lambda}\otimes V_K^+(i(\lambda)),$$
where $V_K^+(a)$ denotes the irreducible $V_K^+$-module labeled by $a\in A$.
\begin{prop}\label{v}
The $V_N\otimes V_K^+$-module $V$ has a unique simple VOA structure extending the 
VOA $V_N\otimes V_K^+$.
\end{prop}
\Pf
The isomorphism $i$ defines the isotropic subspace 
$$C=\{(\lambda,i(\lambda))\mid \lambda \in N'/N\}\subset (N'/N,q_N)\oplus (A,q_A).$$
It follows from Theorem~\ref{abelian-vc} and~\cite{DoLe} that
the direct sum of the irreducible modules of the VOA $V_N\otimes V_K^+$  has the structure
of an abelian intertwining algebra for the finite quadratic space $(N'/N,q_N)\oplus (A,q_A)$.
The proposition follows now from~\cite{Ho-genus}, Theorem~4.3 (or~\cite{DoMa-reductive}).
\eop

\begin{rem}\rm 
The isomorphism type of $V$ could (and in fact does) depend on the chosen 
isomorphism~$i$. The reason is that neither the image of ${\rm Aut}(V_N)$ nor
${\rm Aut}(V_K^+)$ for the induced action on the set of isomorphism classes of irreducible
modules is the full orthogonal group of the corresponding finite quadratic space. This
follows from the observation that in both cases the six orbits of the orthogonal group split
into $16$ respectively $7$ orbits. There are up to automorphism six possibilities for $i$.
They correspond to the VOAs with affine Kac-Moody subVOA
$B_{12,2}$, $B_{6,2}^2$, $B_{4,2}^3$, $B_{3,2}^4$, $B_{2,2}^6$ or $A_{1,4}^{12}$
in Schellekens' list~\cite{ANS} of self-dual VOA candidates of central charge~$24$.


The genus $I\!I_{12,0}(2_{I\!I}^{-10} 4_{I\!I}^{-2})$ of $N$ consists of the two classes 
$\sqrt{2}D_{12}$ and $\sqrt{2}(E_8\oplus D_4)$ which have isomorphic discriminant forms.
If we replace the lattice $N$ in the construction of $V$ by the lattice $\sqrt{2}(E_8\oplus D_4)$ 
the resulting VOAs have the affine Kac-Moody subVOA $A_{8,2}F_{4,2}$, $C_{4,2}A_{4,2}^2$ or $D_{4,4}A_{2,2}^4$.
\end{rem}

We extend the action of the torus $T$ from section~\ref{vn} on $V_N$ and its modules
to $V$ by taking the trivial $T$-action on $V_K^+$ and its modules. 
Note that the $T$-action on $V$ is compatible with the Virasoro module structure.

For the next theorem, we define $f_n=g_n/\eta^{12}$. Explicitly, one has
\begin{equation}\label{function-f}
\begin{array}{lcl}
 f_1 & = & q^{-1} + 12 + 300\,q + 5792\,q^2 + 84186\,q^3 + \cdots,  \\
 f_2 & = & 12 + 288\,q + 5792\,q^2 + 84096\,q^3 + \cdots,           \\
 f_3 & = & 4  + 224\,q + 5344\,q^2 + 81792\,q^3 + \cdots,           \\
 f_4 & = & q^{-1/2} + 44\,q^{1/2} + 1242\,q^{3/2} + 22216\,q^{5/2} +\cdots \\
 f_5 & = & 32\,q^{1/2} + 1152\,q^{3/2} + 21696\,q^{5/2} +\cdots,     \\
 f_6 & = & 12\,q^{1/4} + 520\,q^{5/4} + 10908\,q^{9/4} + \cdots,       \\
 f_7 & = & q^{-1/4} + 90\,q^{3/4} + 2535\,q^{7/4} + 42614\,q^{11/4} +\cdots.
\end{array}
\end{equation}

For $\gamma\in N'/N$, we let $f_\gamma=f_n$ if 
$i(\gamma)$ belongs to the ${\rm Aut}(V_K^+)$-orbit no.~$n$ in Table~\ref{VKP-modules}.


The expression for the $T$-equivariant graded character of $V$ at which we arrive is described
in the following theorem.
\begin{thm}\label{character-lattice}
$$ \chi_V(\tau,{\bf z}) = 
  \sum_{\gamma\in N'/N} f_{\gamma}(\tau)\, \theta_{\gamma}(\tau,{\bf z}). $$
\end{thm}

\Pf 
The 
statement follows from Lemma~\ref{gitterchar} together with the definition
of $V$ and the $f_{\gamma}$.  \eop


\section{The \gkm $\g$}\label{fbmla}

In this section, we construct a new \gkm $\g$ from $V$. 
We determine its simple roots using the singular theta correspondence.
 
\medskip 

Let $\II$ be the even unimodular Lorentzian lattice of rank $2$
and $V_{\II}$ the associated lattice vertex algebra.
Let $V$ be the VOA of the last section. Then the tensor product 
$W=V\otimes V_{\II}$ is a vertex algebra of central charge~$26$.

Let $L=N\oplus \II$. Since $\II$ is unimodular this decomposition gives an isomorphism 
between the discriminant form of $L$ and that of $N$.

\begin{lem}
The isomorphism type of the vertex algebra $W$ does not depend on the isomorphism $i$
used in the definition of $V$.
\end{lem}

\Pf 
{}From the isomorphism $i:(N'/N,q_N) \to (A,-q_A)$ we obtain an isomorphism $i':(L'/L,q_L) \to (A,-q_A)$ 
and $W$ has as $V_L\otimes V_K^+$-module the decomposition 
$$V\ =\ \bigoplus_{\gamma \in  L'/L}   V_\gamma \otimes V_K^+(i'(\gamma)) .$$
Since ${\rm Aut}(L)$ maps surjectively onto the automorphism group of $(L'/L,q_L)$
(Theorem~1.14.2, \cite{Ni-genus}) 
the same holds for the induced action of ${\rm Aut}(V_L)$ on the set of isomorphism types 
of $V_L$-modules. Hence the result of the gluing 
depends up to an automorphism of $V_L$ not on the chosen isomorphism $i'$.
\eop

We remark that if in the construction of $V$ the lattice $N$ is replaced by 
$\sqrt{2}(E_8\oplus D_4)$ then the resulting vertex algebra 
is isomorphic to $W$ because $N\oplus \II \cong \sqrt{2}(E_8\oplus D_4) \oplus \II$.

\medskip
There is an action of the BRST-operator on the tensor product 
of the vertex algebra $W$ of central charge~$26$
with the bosonic ghost vertex superalgebra
$V_{\rm ghost}$ of central charge~$-26$,  
which defines the BRST-cohomology groups  $H^n_{\rm BRST}(W)$. 
The degree one cohomology group $H^1_{\rm BRST}(W)$ has additionally 
the structure of a Lie algebra \cite{FGZ,LZ-gerstenhaber}. 

We can assume that $V$ is defined over the field of real numbers. 
The same holds for the vertex algebra  $\VII$, 
for $V_{\rm ghost}$ and hence for $W$.
\begin{de}\rm
We define the Lie algebra $\g$ as  $H^1_{\rm BRST}(W).$
\end{de}

Then the no-ghost theorem implies 
the following (cf.\ Prop.~3.2, \cite{HS-fakebaby}):

\begin{prop}\label{rootspace}
The Lie algebra $\g$ is a \gkm graded by the lattice $N'\oplus \II=L'$.
Its components $\g(\a)$, for $\a=(\s,r)\in N'\oplus \II$ are
isomorphic to $V_{1-r^2/2}(2\,\s)$ for $\a\not = 0$ and to 
$V_1(0)\oplus {\bf R}^{1,1}\cong {\bf R}^{13,1}$ for $\a=0$.
\phantom{xxxxxxxxxxxxxxx}\hfill \eop
\end{prop} 
The subspace $\g(0)$ of degree $0 \in L'$ is a Cartan subalgebra for $\g$.

We denote the Fourier coefficient of $f_{\gamma}$ at $q^n$ by 
$[f_{\gamma}](n)$.

\begin{thm}\label{rootmult}
For a nonzero vector $\a\in L'$
the dimension of $\g(\a)$ is given by
$$ \dim \g(\a)= [f_{\alpha+L}](-\alpha^2/2) .   $$
The dimension of the Cartan subalgebra is~$14$.
\end{thm}

\Pf Theorem~\ref{character-lattice} and Proposition~\ref{rootspace}.
\eop 

{}From the Fourier expansion of the $f_{\gamma}$ we can read off the real 
roots of $\g$. Recall that we use the isomorphism $i'$ to identify $(L'/L,q_L)$ with $(A,-q_A)$.

\begin{cor}
The real roots of $\g$ are the vectors\\
\hspace*{4mm} $\al \in L$   with $\al^2 = 2$, \\
\hspace*{4mm} $\al \in L'$  with $\al^2 = 1$ and ${i'}(\alpha + L)$ 
belongs to the orbit no.~4 in Table~\ref{VKP-modules}, \\
\hspace*{4mm} $\al \in L'$  with $\al^2 = 1/2$.\\
They all have multiplicity~$1$.
\end{cor}

The reflections in the real roots generate the Weyl group $W$ of $\g$. 

The Weyl group $W$ has a Weyl vector, i.e.~there is a vector $\rho$ in $L' \otimes \R$ 
such that a set of simple roots of $W$ are the roots $\al$ of $W$ satisfying $(\rho, \al) = - \al^2/2$.
The vector $2\rho$ is a primitive norm $0$ vector in $L'$ and $2\rho$ is in $2L' \! \mod L$
(cf.\ Th.~12.1 and~10.4 in \cite{Bo-theta}).

\begin{prop}
The simple roots of the reflection group $W$ form a set of real simple roots for $\g$.
\end{prop}

\Pf
This follows from Cor.~2.4 in \cite{Bo-gen1}.
\eop 

Since $L'$ is Lorentzian $L'\otimes \R$ has two cones of vectors of norm $\leq 0$. The inner product 
of two nonzero vectors in one of the cones is at most $0$ and equal to~$0$ if and 
only if both are positive multiples of the same norm~$0$ vector.  

\begin{prop}
The vectors $2n\rho$, where $n$ is a positive integer, are imaginary simple roots 
of multiplicity $12$.
\end{prop}

\Pf  Since $\rho$ has negative inner product with all real simple roots, $\rho$ lies 
in the fundamental Weyl chamber $C$. We can choose imaginary simple roots lying 
in $C$ (Prop.\ 2.1 in~\cite{Bo-gen1}). It follows that $\rho$ has inner product 
$\leq 0$ with all simple roots. Now write $2n\rho$ as a sum of simple roots with 
positive integral coefficients, i.e.\ $2n\rho = \sum c_i \al_i$. 
Then $0 = \sum c_i (\al_i, 2n\rho ) \leq 0$
so that $(\al_i, 2n\rho ) = 0$ for all $i$. Since $2\rho$ is primitive in $L'$ it follows 
that $\al_i$ is a positive multiple of $2\rho$. Finally all positive multiples of $2\rho$ are
simple roots because the support of a root is connected. The $L$-cosets of the $2n\rho$ 
are mapped by $i'$ to the orbits nos.~1 and~2 in Table~\ref{VKP-modules}.
The constant coefficient of $f_1$ and $f_2$ is $12$ so that the $2n\rho$ all have multiplicity $12$.
\eop

We will see that we have already found a complete set of simple roots for $\g$.

\begin{prop}
The function $F=\sum_{\gamma \in L'/L} f_{\gamma} e^{\gamma}$ is a vector valued modular form 
of weight $-6$ for the Weil representation of ${\rm SL}_2(\Z)$ associated to $L'/L$.
\end{prop}

This follows in principle from the theory of VOAs since the $f_{\gamma}$ are up to 
the factor $1/\eta^{12}$ the characters of the irreducible $V_K^+$-modules
and the VOA $V_K^+$ has a modular tensor category associated
to the finite quadratic space $(A,q_A)$. However, we will give a direct proof.

\Pf 
Since we identify $(L'/L,q_L)$ with $(A,-q_A)$ by $i'$
we have to show that $F =\sum_{a\in A}f_a e^a$, 
where $f_a=f_n$ if $a$ belongs to the ${\rm Aut}(V_K^+)$-orbit 
no.~$n$ in Table~\ref{VKP-modules}, is a vector valued modular form of weight~$-6$ 
with respect to the dual Weil representation of ${\rm SL}_2(\Z)$ for the quadratic space $(A,q_A)$.

The theta function $\Theta_K = \sum_{\mu \in K'/K} \theta_{\mu }e^{\mu}$
transforms under the dual Weil representation of $K'/K$. 
Hence $\Theta_K/(2\Delta)$ where $\Delta = \eta^{24}$ 
is a modular form of weight $-6$ for the dual Weil representation associated to $K'/K$.

Let $H=\{[1]^+,\,[1]^-\}$ be the order~$2$ subgroup of $A$ corresponding to the two $V_K^+$-modules $[1]^+$ and $[1]^-$. 
Then the orthogonal complement $H^{\perp}$ of $H$ in $A$ consists of the set 
of untwisted $V_K^+$-modules denoted by $[n]^{\pm}$, $n=1$,~$\ldots$,~$7$, 
in Table~\ref{VKP-modules} and the quotient $H^{\perp}/H$ is naturally isomorphic to $K'/K$. 

Let $F_K = \sum_{a \in A} F_{K,a} e^a$ be the function with components 
$F_{K,a} =  {\theta}_{\lambda+K}/(2\Delta)$ 
if $a \in H^{\perp}$ is mapped to $\lambda+K$ in  $H^{\perp}/H\cong K'/K$ and $F_{K,a} = 0$ otherwise. 
It follows that $F_K$ is a modular form of weight $-6$ for the dual Weil representation of $(A,q_A)$. 

Let $h(\tau) = 1/\eta(2\tau)^{12}$ and denote by $F_{h/2,0}$ and $F_{-h/4,H}$ 
the lifts of $h/2$ and $-h/4$ on the isotropic subgroups~$0$ and~$H$,
respectively (cf.~\cite{Scheit-Weil}). 
The liftings $F_{h/2,0}$ and $F_{-h/4,H}$
are also modular forms of weight $-6$ for the dual Weil representation of $(A,q_A)$.

Explicit calculations using the equations~(\ref{characters}) together with the identities arising 
from the induced action of ${\rm Aut}(V_K^+)$ on $A$
show that
\[ F = F_K + F_{h/2,0} + F_{-h/4,H}.    \]
\eop

The next result is a consequence of the singular theta correspondence.

\begin{thm}
A set of simple roots for $\g$ is the following: The real simple roots are the vectors\\
\hspace*{4mm} $\al \in L$   with $\al^2 = 2$, \\
\hspace*{4mm} $\al \in L'$  with $\al^2 = 1$ and ${i'}(\alpha + L)$ 
belongs to the orbit no.~4 in Table~\ref{VKP-modules}, \\
\hspace*{4mm} $\al \in L'$ with $\al^2 = 1/2$\\
and which satisfy $(\rho,\al) = -\al^2/2$. 
The imaginary simple roots are the positive integral multiples of $2\rho$ each with multiplicity~$12$.
\end{thm}

\Pf Let $M = L\oplus  I\!I_{1,1} = N \oplus I\!I_{1,1} \oplus I\!I_{1,1}$. 
Then $M'/M$ is isomorphic to $N'/N$ and hence $F$ defines a vector valued modular form for 
the Weil representation of ${M'/M}$. 
The singular theta correspondence associates to $F$ an automorphic product 
$\Psi$ on the Grassmannian of two-dimensional negative definite subspaces of $M\otimes \R$. 
The level one expansion of $\Psi$ is given by 
\[ e((\rho,Z)) \prod_{\al \in L'^+}  \big(1 - e((\al,Z)) \big)^{[f_{\al + L}](-\al^2/2)}  \, . \]
The automorphic form $\Psi$ has singular weight so that the Fourier expansion is 
supported only
on norm~$0$ vectors. Furthermore, $\Psi$ is antisymmetric under the Weyl group $W$. 
It follows that $\Psi$ has the sum expansion 
\[ \sum_{w \in W} \det(w) \, e((w\rho,Z))\prod_{n>0}\big(1-e((2nw\rho,Z))\big)^{12} \, . \]
Now let $\k$ be the generalized Kac-Moody algebra with simple roots as stated in the theorem 
and Cartan subalgebra $L'\otimes \R$. Then the above argument shows that the denominator identity 
of $\k$ is given by 
\[ e^{\rho} \prod_{\al \in L'^+} (1 - e^{\al})^{[f_{\al + L}](-\al^2/2)}  
 = \sum_{w \in W} \det(w) \, w\Big( e^{\rho} \prod_{n>0}(1-e^{2n\rho})^{12} \Big) \, .  \]
Hence $\g$ and $\k$ have the same root multiplicities. 
When we have fixed a Cartan subalgebra and a fundamental Weyl chamber the 
root multiplicities of a generalized Kac-Moody algebra determine the simple roots 
because of the denominator identity. 
It follows that $\g$ and $\k$ have the same simple roots and therefore are isomorphic. 
\eop

\begin{cor} \label{di}
The denominator identity of $\g$ is 
\[ e^{\rho} \prod_{\al \in L'^+} (1 - e^{\al})^{[f_{\al + L}](-\al^2/2)}  
 = \sum_{w \in W} \det(w) \, w\Big( e^{\rho} \prod_{n>0}(1-e^{2n\rho})^{12} \Big) \, .  \]
\eop
\end{cor}

\medskip

We finish with two remarks.
\begin{rem}\rm
The Lie algebra $\g$ can also be constructed by orbifolding the fake monster algebra. 
This is a generalized Kac-Moody algebra describing the physical states of a bosonic string 
moving on a 26-dimensional torus. An extension of $\mathrm{Co}_0$ acts on this Lie algebra. 
By taking the trace of an element over the identity $\Lambda^*(E)=H^*(E)$ we obtain 
a twisted denominator identity. 
A suitable lift of an element of class $2C$ in $\mathrm{Co}_0$ gives the identity in 
Corollary \ref{di}. 
For more details see \cite{Bo-lie}, Section 13 and \cite{Scheithauer-classification}, Section 10.
\end{rem}

\begin{rem}\rm
The above method for the construction of $\g$ can also be used to construct 
the fake baby monster algebra~\cite{Bo-lie,HS-fakebaby}.
In this case one takes for $K$ the rank~$8$ lattice $\sqrt{2}E_8$.
The VOA $V_K^+$ has an abelian intertwining algebra based on a finite quadratic space
$(A,q_A)$ with $A$ a $2$-elementary group of order $2^{10}$. The automorphism group $O^+(10,2)$
of $V_K^+$ equals in this case the isomorphism group of $(A,q_A)$~\cite{Griess,Shimakura1}.
For the lattice $N$, one can take any of the $17$ lattices in the corresponding genus 
$I\!I_{16,0}(2^{+10}_{I\!I})$.
The resulting VOAs $V$ which here clearly do not depend on the chosen isomorphisms~$i$ 
are the $17$ VOAs occurring in Schellekens' list of self-dual VOAs of central
charge~$24$ having a Lie algebra $V_1$ of rank~$16$~\cite{ANS}.
The resulting vertex algebras $W$ and the corresponding Lie algebras are again isomorphic.
In the paper~\cite{HS-fakebaby}, we started with $V$ belonging to the
affine Kac-Moody VOA $A_{1,2}^{16}$.
\end{rem}


\small


\providecommand{\bysame}{\leavevmode\hbox to3em{\hrulefill}\thinspace}
\providecommand{\MR}{\relax\ifhmode\unskip\space\fi MR }
\providecommand{\MRhref}[2]{%
  \href{http://www.ams.org/mathscinet-getitem?mr=#1}{#2}
}
\providecommand{\href}[2]{#2}

\end{document}